\documentclass[a4paper,leqno,12pt,twoside]{amsart}
\usepackage{amsmath,amscd,amssymb,upgreek,ifthen} 
\usepackage[all]{xy}
\usepackage{graphicx}
\usepackage{yhmath}
\usepackage[french]{babel}

\usepackage[applemac] {inputenc}

\def\C{\mathbb{C}}
\def\H{\mathbb{H}}
\def\N{\mathbb{N}}
\def\U{\mathbb{U}}
\def\Q{\mathbb{Q}}
\def\Z{\mathbb{Z}}

\def\D{\mathbb{D}}
\def\A{\mathbb{A}}
\def\S{\mathbb{S}}
\def\Ann{\operatorname{Ann}}
\def\Hom{\operatorname{Hom}}
\def\Aut{\operatorname{Aut}}

\def\pgcd{\operatorname{pgcd}}

\newcommand{\Exemple}[1]{\begin{Exple}\emph{#1}\end{Exple}}
\newcommand{\Remarque}[1]{\begin{Rem}\emph{#1}\end{Rem}}
\newcommand{\Theoreme}[1]{\begin{The}#1\end{The}}
\newcommand{\Proposition}[1]{\begin{Prop}#1\end{Prop}}
\newcommand{\Lemme}[1]{\begin{Lem}#1\end{Lem}}
\newcommand{\Corollaire}[1]{\begin{Cor}#1\end{Cor}}
\newcommand{\Definition}[1]{\begin{Def}#1\end{Def}}
\newcommand{\Preuve}[1]{\noindent \textbf{Preuve : }#1 \hfill$\blacksquare$\\}
 
\newcommand{\application}[5]{
\ifthenelse{\equal{#1}{0}}{}{#1:}\begin{array}[t]{ccl}
	#2 & \longrightarrow & #3
	\ifthenelse{\equal{#4}{0}}{}{ \\ #4 & \longmapsto & #5} 
\end{array}}

\begin{document}
\title{SURFACES DE STEIN ASSOCIÉES AUX SURFACES DE KATO INTERMÉDIAIRES}
\author{Laurent BATTISTI}  
\address{{\it Laurent BATTISTI}:
LATP-UMR(CNRS) 6632, 
CMI-Universit\'e d'Aix-Marseille I, 39, rue Joliot-Curie, F-13453 Marseille Cedex 13, France.}
\email{battisti@cmi.univ-mrs.fr}
\thanks{Le financement de cette recherche est assuré par la Région Provence-Alpes-Côte d'Azur dans le cadre d'une bourse doctorale régionale.}

\newtheorem{The}{Théorème}[section]
\newtheorem{Prop}[The]{Proposition}
\newtheorem{Def}[The]{Définition}
\newtheorem{Lem}[The]{Lemme}
\newtheorem{Cor}[The]{Corollaire}
\newtheorem{Rem}[The]{Remarque}
\newtheorem{Exple}[The]{Exemple}

\newcommand{\Chapitre}[1]{\section{#1}}

\newcommand{\SousChapitre}[1]{\subsection{#1}}
\newcommand{\SousSousChapitre}[1]{\subsubsection{#1}}
\newcommand{\Paragraphe}[1]{\paragraph{#1}}
\newcommand{\SousParagraphe}[1]{\subparagraph{#1}}

\maketitle

\begin{abstract}
Let $S$ be an intermediate Kato surface, $D$ the divisor consisting of all rational curves of $S$, $\widetilde{S}$ the universal covering of $S$ and $\widetilde{D}$ the preimage of $D$ in $\widetilde{S}$. We prove two results about the surface $\widetilde{S}\setminus \widetilde{D}$: it is Stein (which was already known when $S$ is either a Enoki or a Inoue-Hirzebruch surface) and we give a necessary and sufficient condition so that its holomorphic tangent bundle is holomorphically trivialisable.

-----

Soient $S$ une surface de Kato interm\'ediaire, $D$ le diviseur form\'e des courbes rationnelles de $S$, $\widetilde{S}$ le rev\^etement universel de $S$ et $\widetilde{D}$ la pr\'eimage de $D$ dans $\widetilde{S}$. On donne deux r\'esultats concernant la surface $\widetilde{S}\setminus \widetilde{D}$, \`a savoir qu'elle est de Stein (ce qui \'etait connu dans le cas o\`u $S$ est une surface d'Enoki ou d'Inoue-Hirzebruch) et on donne une condition nécessaire et suffisante pour que son fibré tangent holomorphe soit holomorphiquement trivialisable. 
\end{abstract}

\Chapitre{Introduction}
Les surfaces de la classe $\rm{VII}$ de Kodaira sont les surfaces complexes compactes dont le premier nombre de Betti vaut $1$ ; on appelle surface de la classe $\rm{VII}_{0}$ une surface de la classe $\rm{VII}$ qui est minimale. Le cas de ces surfaces dont le second nombre de Betti $b_{2}$ est nul est entièrement compris, il s'agit nécessairement d'une surface de Hopf ou d'une surface d'Inoue et le cas $b_{2}>0$ est toujours étudié actuellement ; il a été conjecturé qu'elles contiennent toutes une coquille sphérique globale. La preuve de ce résultat terminerait la classification des surfaces complexes compactes.

Les surfaces à coquille sphérique globale, qui nous intéressent ici, peuvent être obtenues selon un procédé dû à  Kato (voir \cite{Kato:1978aa}), que l'on rappelle dans la section suivante. Ces surfaces se divisent en trois classes, les surfaces d'Enoki, d'Inoue-Hirzebruch et enfin les surfaces intermédiaires.

Étant donnés une surface minimale $S$ à coquille sphérique globale, $D$ le diviseur maximal de $S$ formé des $b_{2}(S)$ courbes rationnelles de $S$ et $\varpi : \widetilde{S}\rightarrow S$ le revêtement universel de $S$, nous allons démontrer que $\widetilde{S}\setminus\widetilde{D}$ (où $\widetilde{D}=\varpi^{-1}(D)$) est une variété de Stein. Ce résultat était déjà connu pour les surfaces d'Enoki et d'Inoue-Hirzebruch ; nous allons le montrer dans le cas des surfaces intermédiaires. Dans la dernière partie et toujours dans le cas des surfaces intermédiaires, on donne une condition pour que le fibré tangent holomorphe de la variété $\widetilde{S}\setminus\widetilde{D}$ soit holomorphiquement trivialisable, à savoir que la surface $S$ soit d'indice $1$.\\

\Chapitre{Préliminaires}
On dit qu'une surface compacte $S$ contient une coquille sphérique globale s'il existe une application qui envoie biholomorphiquement un voisinage de la sphère $\S^{3} \subset \C^{2}\setminus \{0\}$ dans $S$ et telle que le complémentaire dans $S$ de l'image de la sphère par cette application soit connexe. \\
Toute surface contenant une coquille sphérique globale peut être obtenue de la façon suivante : étant données une succession finie d'éclatements $\pi_{1}, ..., \pi_{n}$ de la boule unité $B$ de $\C^{2}$ au-dessus de $0$ et  ${\pi:=\pi_{1}\circ \cdots \circ \pi_{n} : B^{\pi} \rightarrow B}$ la composée de ces éclatements, ainsi qu'une application ${\sigma : \overline{B} \rightarrow B^{\pi}}$ biholomorphe sur un voisinage de $\overline{B}$, on recolle les deux bords de $\Ann(\pi,\sigma):=B^{\pi}\setminus \sigma(\overline{B})$ à l'aide de l'application $\sigma \circ \pi$ : 
\begin{figure}[!h] 
   \centering
   \includegraphics[scale=0.57]{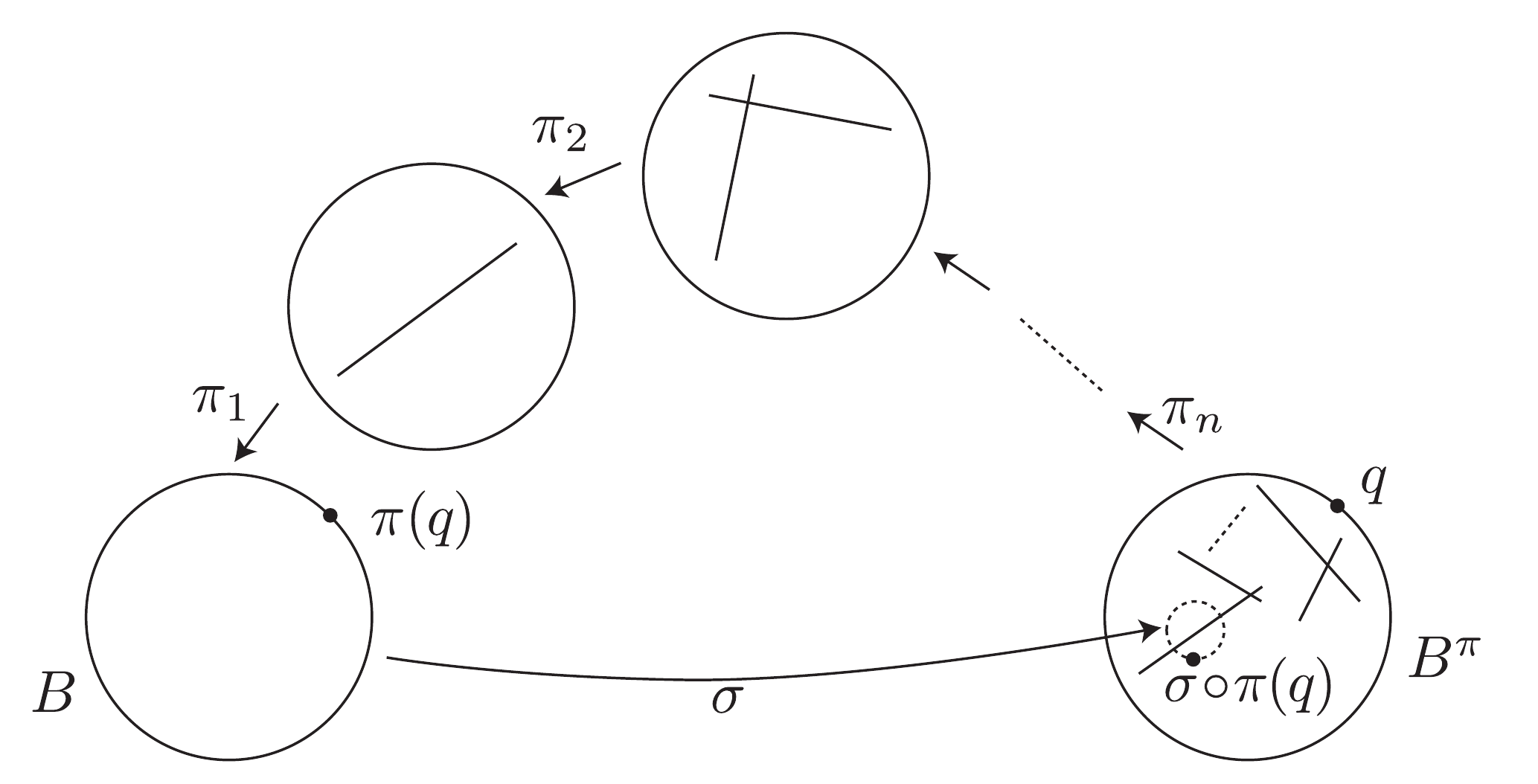}\\ 
\end{figure}

La surface obtenue possède un groupe fondamental isomorphe à $\Z$ et son second nombre de Betti est égal à $n$ (voir \cite{Dloussky:1984aa}). Il s'agit d'une construction due à Kato \cite{Kato:1978aa}. Dans la suite, on appellera surface de Kato une surface complexe compacte minimale contenant une coquille sphérique globale, dont le second nombre de Betti est non nul.

Dans \cite{Dloussky:1984aa}, Dloussky étudie le germe contractant d'application holomorphe $\varphi = \pi\circ \sigma : B \rightarrow B$ associé à la construction précédente. Ce germe détermine à isomorphisme près la surface étudiée (proposition 3.16 loc. cit.).  \\

Soit $S$ une surface de Kato ; on note $D$ le diviseur maximal de $S$ formé des $b_{2}(S)$ courbes de $S$, $\widetilde{S}$ le revêtement universel de $S$ et $\widetilde{D}$ la préimage de $D$ dans $\widetilde{S}$.

Suivant les notations de \cite{Dloussky:1984aa}, on obtient la surface $\widetilde{S}$ en recollant une infinité d'anneaux $A_{i}$ ($i \in \Z$) isomorphes à $\Ann(\pi,\sigma)$, en identifiant le bord pseudo-concave de $A_{i}$ au bord pseudo-convexe de $A_{i+1}$ via l'application $\sigma \circ \pi$. La surface $\widetilde{S}$ possède deux bouts, notés $\underline{0}$ et $\infty$, le bout $\underline{0}$ possédant une base de voisinages ouverts strictement pseudo-convexes (les $\bigcup_{i \geqslant j}A_{i}$ pour $j \in \Z$) et le second une base de voisinages strictement pseudo-concaves (les $\bigcup_{i \leqslant j}A_{i}$ pour $j \in \Z$). Enfin on définit un automorphisme $G$ de $\widetilde{S}$ en posant $G(z_{i}):=z_{i+1}$ où $z_{i}$ et $z_{i+1}$ sont les images dans $A_{i}$ et $A_{i+1}$ respectivement d'un même point $z\in \Ann(\pi,\sigma)$.
\begin{figure}[!h] 
   \centering
   \includegraphics[scale=0.57]{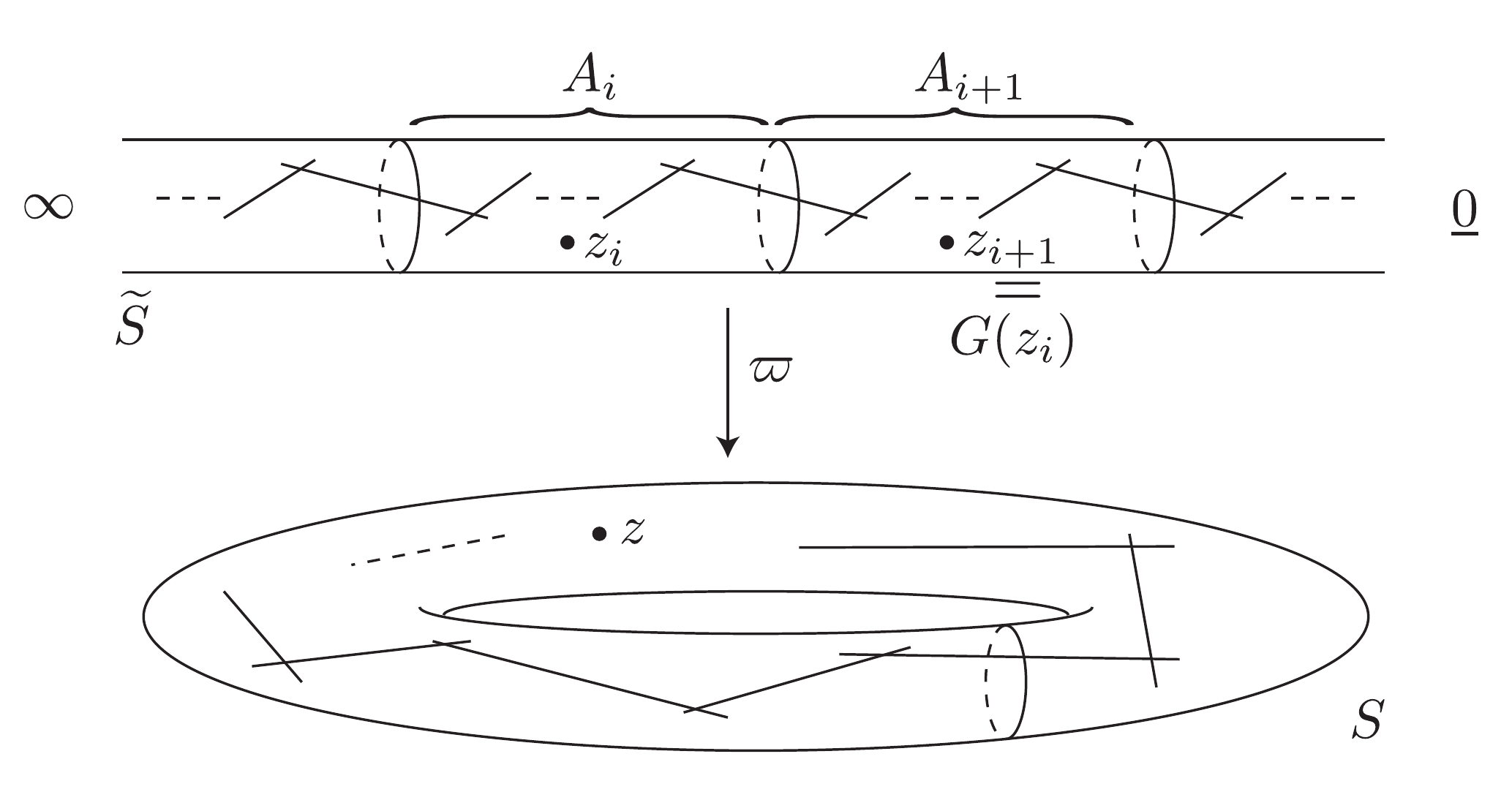}\\ 
\end{figure}

Fixons une courbe compacte $C$ de $\widetilde{S}$ avec $C \subset A_{0}$. On note $(\widehat{S}_{C},p_{C})$ l'effondrement de $\widetilde{S}$ sur la courbe $C$, c'est-à-dire la donnée d'une surface $\widehat{S}_{C}$ n'ayant qu'un bout, d'une application holomorphe $p_{C}$ de $\widetilde{S}$ dans $\widehat{S}_{C}$, biholomorphe sur un voisinage du bout $\infty$ dans $\widetilde{S}$ sur un voisinage du bout de $\widehat{S}_{C}$, telles que $\widehat{C}=p_{C}(C)$ soit une courbe d'auto-intersection $-1$.

La proposition 3.4 de \cite{Dloussky:1984aa} nous assure l'existence d'une telle application $p_{C}$ pour toute courbe compacte $C$ de $\widetilde{S}$, et d'un point $\widehat{0}_{C} \in \widehat{C}$ tel que $p_{C}$ soit également biholomorphe entre $\widetilde{S}\setminus p_{C}^{-1}(\widehat{0}_{C})$ et $\widehat{S}_{C}\setminus \{\widehat{0}_{C}\}$. 

De plus, la restriction de $p_{C}$ au complémentaire de $\widetilde{D}$ est un biholomorphisme entre $\widetilde{S}\setminus \widetilde{D}$ et $\widehat{S}_{C} \setminus p_{C}(\widetilde{D})$. Enfin, il existe une application holomorphe $F_{C}$ de $\widehat{S}_{C}\setminus \{\widehat{0}_{C}\}$ dans lui-même, contractante en $\widehat{0}_{C}$, conjuguée à $\varphi$ et biholomorphe sur $\widehat{S}_{C}\setminus p_{C}(\widetilde{D})$. 
~\\

\Chapitre{La variété $\widetilde{S}\setminus \widetilde{D}$ est de Stein}

Les surfaces de Kato se divisent en trois classes : les surfaces d'Enoki, d'Inoue-Hirzebruch et enfin les surfaces intermédiaires (voir \cite{Dloussky:2003fk}).

Dans le cas des surfaces d'Inoue-Hirzebruch et celles d'Enoki, le fait que $\widetilde{S}\setminus\widetilde{D}$ soit de Stein est déjà connu : pour une surface d'Inoue-Hirzebruch, la variété $\widetilde{S}\setminus\widetilde{D}$ est un domaine de Reinhardt holomorphiquement convexe (voir \cite{Zaffran:2001aa}, proposition 2.2) tandis que pour une surface d'Enoki, on a $\widetilde{S}\setminus\widetilde{D}\cong \C^{*}\times\C$ qui sont bien dans chaque cas des variétés de Stein. Il reste donc à étudier le cas des surfaces intermédiaires.
~\\

Favre a donné dans \cite{Favre:2000uq} des formes normales pour les germes contractants d'applications holomorphes et on peut en particulier donner la forme du germe associé à une surface intermédiaire, à savoir qu'une telle surface est associée au germe $\varphi$ de $(\C^{2},0) \rightarrow (\C^{2},0)$ donné par \begin{equation}\label{formeNormale}(z,\zeta) \mapsto (\lambda \zeta^sz+P(\zeta)+c_{\frac{sk}{k-1}}\zeta^{\frac{sk}{k-1}},\zeta^{k})\end{equation} où $\lambda \in \C^{*}$, $k, s \in \N$ avec $k>1$ et $s>0$,  et $P(\zeta)=c_{j}\zeta^{j}+...+c_{s}\zeta^{s}$ avec les conditions suivantes : $0<j<k$, $ j\leqslant s$, $c_{j}=1$, $c_{\frac{sk}{k-1}}=0$ quand $\frac{sk}{k-1} \not \in \Z$ ou $\lambda \neq 1$ et enfin $\pgcd\{k,m~|~c_{m}\neq 0\} =1$. On trouve dans \cite{Oeljeklaus:2008aa} une condition pour que deux tels germes soient conjugués (et déterminent donc deux surfaces isomorphes).\\

L'objectif de cette section est de démontrer, dans le cas de surfaces intermédiaires, le \Theoreme{La surface $\widetilde{S}\setminus \widetilde{D}$ est de Stein.}

Dans un premier temps (section \ref{redProb}), on montre qu'il est suffisant de se ramener à la situation du théorème \ref{unionStein} énoncé ci-dessous. Pour cela, nous allons écrire notre surface comme réunion croissante d'ouverts et nous verrons que seule une hypothèse manque \emph{a priori} pour pouvoir effectivement appliquer ce théorème, à savoir que chaque paire constituée de deux tels ouverts consécutifs est de Runge. C'est dans la section \ref{sectionSuivante} qu'on prouve que cette hypothèse est bien vérifiée.

\subsection{\label{redProb}Réduction du problème.}Reprenons les notations précédentes et donnons-nous un germe de la forme (\ref{formeNormale}). On regarde la surface intermédiaire $S$ associée et on choisit une courbe $C$ de $\widetilde{S}$ donnée par la proposition 3.16 de \cite{Dloussky:1984aa} ; quitte à renuméroter les $A_{i}$ on suppose que $C \subset A_{0}$. Notre objectif est de prouver que la variété $\widehat{S}_{C} \setminus p_{C}(\widetilde{D})$ est de Stein, en utilisant le théorème suivant (voir \cite{Gunning:1965kx}, théorème 10 p.~215) : 

\Theoreme{\label{unionStein}Soient $X$ un espace analytique complexe et $(X_{i})_{i\in\N}$ une suite croissante de sous-espaces de $X$ qui soient de Stein. Supposons que $X = \bigcup X_{i}$ et que chaque paire $(X_{i+1},X_{i})$ est de Runge, i.e. l'ensemble $\mathcal{O}(X_{i})|_{X_{i+1}}$ des restrictions à $X_{i}$ des applications holomorphes sur $X_{i+1}$ est dense dans $\mathcal{O}(X_{i})$. Alors $X$ est de Stein.}

Notons : 
\begin{itemize}
\item[-]$\widehat{A}_{i}:=p_{C}(A_{i})$ pour tout $i\in\Z$ et 
\item[-]$\mathcal{A}_{i}:=p_{C}(\bigcup_{j\geqslant i}A_{j})$ pour $i \leqslant 0$,
\end{itemize}
de sorte qu'on a $\mathcal{A}_{i} \subset \mathcal{A}_{i-1}$ et $\widehat{S}_{C} \setminus p_{C}(\widetilde{D}) = \displaystyle \bigcup_{i\leqslant 0}\mathcal{A}_{i}\setminus p_{C}(\widetilde{D})$. 

Chaque $\mathcal{A}_{i}\setminus p_{C}(\widetilde{D})$ est strictement pseudo-convexe, donc de Stein. De plus, on a $F_{C}(\widehat{A}_{i})=\widehat{A}_{i+1}$ pour $i\leqslant -1$, car le diagramme   \[\xymatrix{
    \widetilde{S} \ar[r]^{\displaystyle G} \ar[d]_{\displaystyle p_{C}} & \widetilde{S} \ar[d]^{\displaystyle p_{C}} \\
    \widehat{S}_{C} \ar[r]_{\displaystyle F_{C}} & \widehat{S}_{C}
  }\]
 est commutatif (c.f. \cite{Dloussky:1984aa}, proposition 3.9). Ainsi, on a \begin{equation}\label{germIsom}F_{C}(\mathcal{A}_{i-1}\setminus p_{C}(\widetilde{D}))=\mathcal{A}_{i}\setminus p_{C}(\widetilde{D}) 
\end{equation}

Supposons établi le fait que la paire $(\mathcal{A}_{0}\setminus p_{C}(\widetilde{D}), F_{C}(\mathcal{A}_{0}\setminus p_{C}(\widetilde{D})))$ est de Runge. Alors la paire $(\mathcal{A}_{-1}\setminus p_{C}(\widetilde{D}),\mathcal{A}_{0}\setminus p_{C}(\widetilde{D}))$ est automatiquement de Runge par l'égalité (\ref{germIsom}) ci-dessus, et par récurrence chaque paire $(\mathcal{A}_{i-1}\setminus p_{C}(\widetilde{D}),\mathcal{A}_{i}\setminus p_{C}(\widetilde{D}))$ est de Runge. Nous sommes alors en mesure d'appliquer le théorème \ref{unionStein} qui nous dit que la réunion des $\mathcal{A}_{i}\setminus p_{C}(\widetilde{D})$ est de Stein.

Le problème est donc ramené à montrer que le couple $(\mathcal{A}_{0}\setminus p_{C}(\widetilde{D}), F_{C}(\mathcal{A}_{0}\setminus p_{C}(\widetilde{D})))$ est de Runge.

\Remarque{L'ensemble $\mathcal{A}_{0} \setminus p_{C}(\widetilde{D})$ est biholomorphe à une boule ouverte centrée en $0$ privée d'une droite complexe. En effet, on peut écrire $\varphi = \pi \circ \sigma$ où $\pi$ est une succession d'éclatements de la boule au-dessus de $0\in\C^{2}$, $\sigma : \overline{B}\rightarrow \pi^{-1}(B)$ est une application définie sur un voisinage de $\overline{B}$ et biholomorphe sur son image, et $\varphi$ est de la forme normale (\ref{formeNormale}). Par le choix de la courbe $C$, la proposition 3.16 p.~33 de \cite{Dloussky:1984aa} nous donne l'isomorphisme $\mathcal{A}_{0} \setminus p_{C}(\widetilde{D})\cong B\setminus\varphi^{-1}(0)$ et en utilisant la forme de $\varphi$, on voit que $\varphi^{-1}(0)=\{\zeta=0\}$. }

Finalement, démontrer que $(\mathcal{A}_{0}\setminus p_{C}(\widetilde{D}),F_{C}(\mathcal{A}_{0}\setminus p_{C}(\widetilde{D})))$ est de Runge revient à prouver que c'est le cas de la paire $\left(B\setminus \{\zeta = 0\}, \varphi(B\setminus \{\zeta = 0\})\right)$ pour une boule $B \subset \C^{2}$ centrée en $0$ (en notant $(z,\zeta)$ les coordonnées de $\C^{2}$). C'est l'objet de la section suivante.

\subsection{\label{sectionSuivante}La paire $\left(B\setminus \{\zeta = 0\}, \varphi(B\setminus \{\zeta = 0\})\right)$ est de Runge}

\noindent Etant donné un germe $\varphi$ de la forme (\ref{formeNormale}),  introduisons en premier lieu quelques notations :\\

\begin{itemize}
\item [1.] Remarquons tout d'abord que chaque point de $\C\times \Delta^{*}$  possède exactement $k$ antécédents par $\varphi$, où $\Delta^{*}$ est le disque unité ouvert de $\C$ privé de $0$. Notons $g$ l'automorphisme de $\C \times \Delta^{*}$ suivant : \[g : (z,\zeta) \mapsto \left(\varepsilon^{-s}z+\displaystyle \frac{P(\zeta)-P(\varepsilon \zeta)}{\lambda \varepsilon^{s}\zeta^{s}},\varepsilon \zeta\right)\] où $\varepsilon$ est une racine primitive $k$-ième de l'unité, de sorte que $\varphi \circ g = \varphi$. Pour tout $\ell \in \Z$, on a \[g^{\ell}(z,\zeta)=\left((\varepsilon^{\ell})^{-s}z+\displaystyle \frac{P(\zeta)-P(\varepsilon^{\ell}\zeta)}{\lambda (\varepsilon^{\ell})^{s}\zeta^{s}},\varepsilon^{\ell}\zeta\right)\] et $g^{\Z} \cong \Z/{k}\Z$. L'automorphisme $g$ permute les antécédents d'un même point de l'application $\varphi$.

\item[2.] On notera également $q(z,\zeta)$ le polynôme $\displaystyle z\prod_{\ell=1}^{k-1}a_{\ell}(z,\zeta)\zeta^{n_{\ell}}$ où $a_{\ell}(z,\zeta)$ est la première composante de $g^{\ell}(z,\zeta)$ et $n_{\ell}=s-\min\{n | c_{n}(1-(\varepsilon^{\ell})^{n}) \neq 0\}$, qui est bien défini et positif ou nul vu la dernière hypothèse sur les coefficients de $P$, à savoir $\pgcd\{k,m~|~c_{m}\neq 0\} =1$. Le polynôme $q(z,\zeta)$ est en particulier de la forme $q(z,\zeta)=z(c+\epsilon(z,\zeta))$ où $\epsilon(z,\zeta) \underset{(z,\zeta) \rightarrow (0,0)}{\xrightarrow{\hspace*{1cm}}} 0$ et $c \neq 0$. \\

\item[3.] Pour $\eta >0$, on note $U_{\eta}$ l'ouvert $\{ (z,\zeta) \in \C^{2} ~|~ |q(z,\zeta)| < \eta \}$. Soient $a$, $b$ et $c$ trois réels strictement positifs, on définit les ensembles \[K_{a,b}:= \{(z,\zeta) \in \C^{2}~|~|z|^{2}+|\zeta|^{2}\leqslant a^{2}, |\zeta| \geqslant b \} = \overline{B(0,a)} \cap \{|\zeta| \geqslant b\} \] et  \[L_{a,b,c}:= \overline{\D(0,a)} \times \overline{\A_{b,c}}\] (où $\A_{b,c}$ est l'anneau ouvert centré en $0$ de rayons $b<c$). Enfin, on pose \[\mathcal{K}_{a,b}:=\bigcup_{\ell=0}^{k-1}g^{\ell}(K_{a,b})\] et \[ \mathcal{L}_{a,b,c}:=\bigcup_{\ell=0}^{k-1}g^{\ell}(L_{a,b,c}).\]
\end{itemize}

\Remarque{\label{rq34}Pour $a,b$ et $c$ assez petits, les compacts $g^{\ell}(L_{a,b,c})$ (resp. $g^{\ell}(K_{a,b})$) sont disjoints deux à deux : ceci est une conséquence du fait que la fonction $\varphi$ est localement injective autour de l'origine de $\C^{2}$, ce qui est démontré, par exemple, dans \cite{Dloussky:2003fk}, section 5. En particulier, les ensembles $\mathcal{L}_{a,b,c}$ et $\mathcal{K}_{a,b}$ possèdent chacun $k$ composantes connexes.\\
D'autre part, on a $K_{a',b'} \subset L_{a,b,c}$ pour $b'\geqslant b$ et $a' \leqslant \min\{a,c\}$, ce qui entraîne notamment $\mathcal{K}_{a',b'} \subset \mathcal{L}_{a,b,c}$.\\
\begin{figure}[!h] 
   \centering
   \includegraphics[scale=0.57]{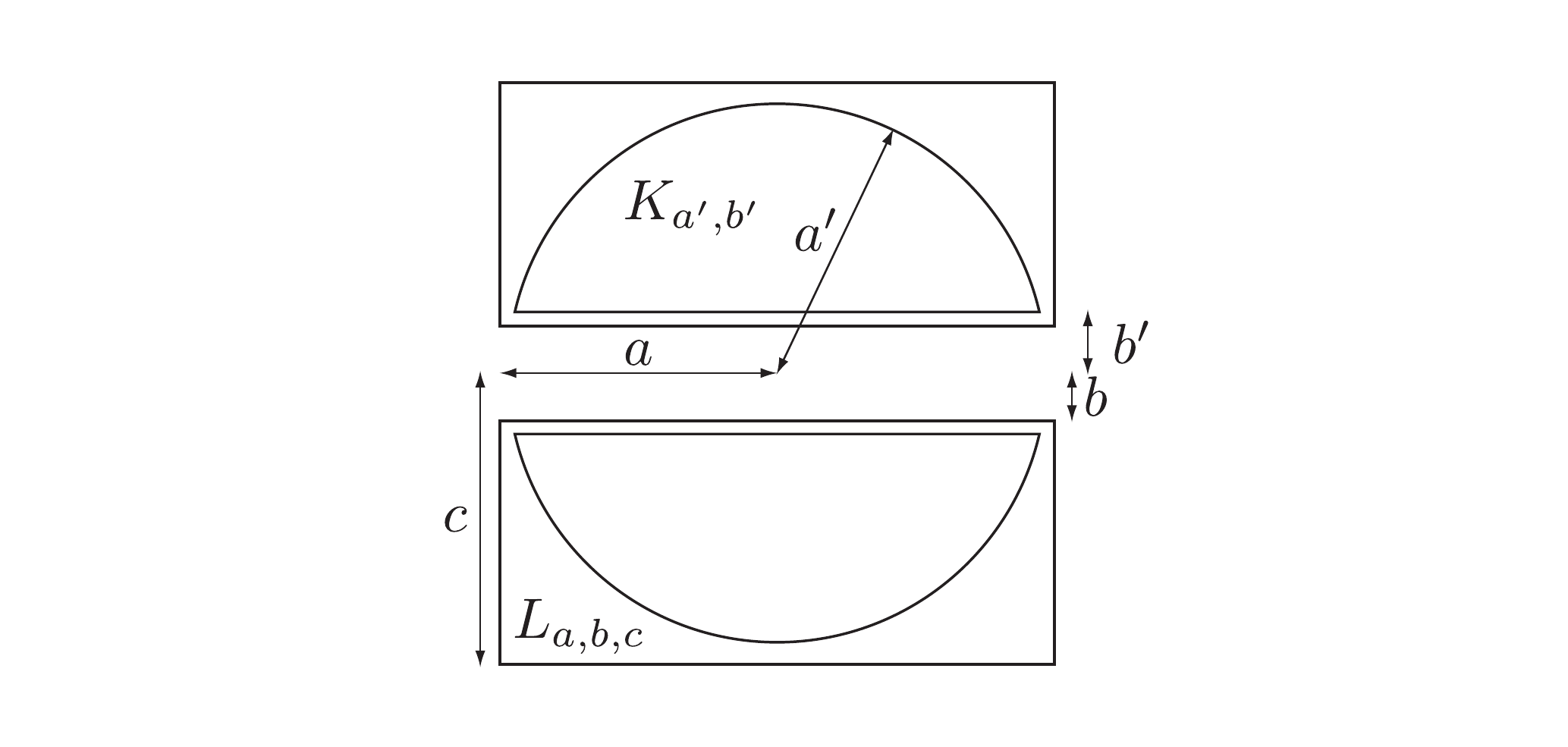}\\ 
\end{figure}\\
Enfin, pour $\eta >0$ fixé, il existe $A_{\eta}>0$ tel que pour tous réels $t$ et $\delta$ avec $0<t<\delta <A_{\eta}$, on ait $\mathcal{L}_{\delta,t,\delta} \subset U_{\eta}$ : en calculant $|q(z,\zeta)|$ pour $(z,\zeta) \in \mathcal{L}_{\delta, t, \delta}$ on voit qu'il suffit de choisir $\delta$ assez petit pour avoir 
\begin{equation}\label{cond1}|\delta |\prod_{\ell=1}^{k-1}(|\delta^{n_{\ell}+1}|+2(|c_{s-n_{\ell}}|+|c_{s-n_{\ell}+1}|\delta+...+|c_{s}|\delta^{n_{\ell}})/\lambda)<\eta.
\end{equation}}

On appelle $V_{\eta,\delta}$ l'ensemble $U_{\eta}\cap\{|\zeta|\leqslant\delta\}$. \\

\Proposition{\label{holConv} Pour $\delta>0$ assez petit et pour tout $\varepsilon_{1} \in ]0,\delta[$, le compact $\mathcal{K}_{\delta,\varepsilon_{1}}$ est holomorphiquement convexe.}
\Preuve{En premier lieu, remarquons que l'enveloppe holomorphiquement convexe de $V_{\eta, \delta}$ est l'adhérence $\overline{V}_{\eta, \delta}$ de cet ensemble. On note :

\begin{itemize}
\item[-] $\widehat{\mathcal{K}}_{\delta,\varepsilon_{1}}$ l'enveloppe holomorphiquement convexe de $\mathcal{K}_{\delta,\varepsilon_{1}}$,
\item[-] $\widehat{\mathcal{K}}^{\ell}_{\delta,\varepsilon_{1}}$ (resp. $\smash{\overline{V}}^{\ell}_{\eta,\delta}$) la composante connexe de $\widehat{\mathcal{K}}_{\delta,\varepsilon_{1}}$ (resp. $\overline{V}_{\eta,\delta}$) qui contient $g^{\ell}(K_{\delta,\varepsilon_{1}})$, pour $\ell \in \{0,...,k-1\}$.\\
\end{itemize}

\noindent \textbf{Étape 1 :} Montrons tout d'abord que pour $\eta$ et $\delta$ assez petits et pour tout $\varepsilon_{1}<\delta$, on a $\smash{\overline{V}}^{0}_{\eta,\delta} \cap \mathcal{K}_{\delta,\varepsilon_{1}}=K_{\delta,\varepsilon_{1}}$, autrement dit que la composante connexe de $\overline{V}_{\eta,\delta}$ qui contient $K_{\delta,\varepsilon_{1}}$ ne rencontre aucune autre composante de $\mathcal{K}_{\delta,\varepsilon_{1}}$.\\

\noindent Soient $\delta>\varepsilon_{1}>0$. Pour $\ell \in \{1,..., k-1\}$, on a \[g^{\ell}(K_{\delta,\varepsilon_{1}}) \subset g^{\ell}(L_{\delta,\varepsilon_{1},\delta})=\{(z,\zeta) \in \C^{2} ~|~ |a_{k-\ell}(z,\zeta)| \leqslant \delta, |\zeta|\in [\varepsilon_{1},\delta] \}.\]
En particulier, pour $(z,\zeta) \in g^{\ell}(L_{\delta,\varepsilon_{1},\delta})$, on a $z=\displaystyle \frac{P(\varepsilon^{k-\ell}\zeta)-P(\zeta)}{\lambda \zeta^{s}}+w$ où $|w|\leqslant \delta$. En développant, cette égalité devient
\[z = \begin{array}[t]{rl}
\lambda^{-1}\zeta^{-n_{\ell}}\left(c_{s-n_{\ell}}((\varepsilon^{-\ell})^{s-n_{\ell}}-1)+\right.&\!\!\!\!\! c_{s-n_{\ell}+1}((\varepsilon^{-\ell})^{s-n_{\ell}+1}-1)\zeta+ ... \\
& \left. ...+c_{s}((\varepsilon^{-\ell})^{s}-1)\zeta^{n_{\ell}}\right)+w.
\end{array}\]
Autrement dit, $z$ est de la forme 
\[\lambda^{-1}\zeta^{-n_{\ell}}((c_{s-n_{\ell}}((\varepsilon^{-\ell})^{s-n_{\ell}}-1)+\zeta R_{\ell}(\zeta,w))\] où $R_{\ell}$ est un polynôme et par définition de $n_{\ell}$, le terme $c_{s-n_{\ell}}((\varepsilon^{-\ell})^{s-n_{\ell}}-1)$ est non nul.\\

Grâce à cette dernière expression de $z$, on voit que lorsque $n_{\ell}>0$, pour n'importe quelle constante $C>0$ et lorsque $\delta$ est assez petit, tout élément $(z,\zeta) \in g^{\ell}(L_{\delta,\varepsilon_{1},\delta})$ vérifie $|z|>C$. 

Dans le cas où $n_{\ell}=0$ (donc $c_{s} \neq 0$), on a $|z| = |\lambda^{-1}(c_{s}((\varepsilon^{-\ell})^{s}-1))+w|$ est supérieur à une constante non nulle pour $\delta$ assez petit.\\

Posons alors $\alpha:=\frac{1}{2}\underset{\ell}{\min}\{ |c_{s} ((\varepsilon^{\ell})^{s}-1)| ~|~ \ell s \not \equiv 0[k] \}$ si $c_{s} \neq 0$ et $\alpha:=1$ sinon. Par ce qui précède, il existe une constante $A>0$ telle que pour tous $\delta<A$ et $\varepsilon_{1}<\delta$ on ait, pour chaque $\ell \in \{1,..., k-1\}$ et tout élément $(z,\zeta)$ de $g^{\ell}(L_{\delta,\varepsilon_{1},\delta})$, l'inégalité 
\begin{equation} \label{cond2}
|z| \geqslant  \alpha. 
\end{equation}

\noindent Fixons désormais $\eta>0$ vérifiant les deux conditions suivantes : 
\begin{itemize} 
\item[1.] $2\eta/|c| < \alpha$ (où $c\neq 0$ est le facteur de $z$ dans le développement limité de $q$ en $(0,0)$, à savoir $q(z,\zeta)=z(c+\epsilon(z,\zeta))$), et 
\item[2.] $|c+\epsilon(z,\zeta)|>|c|/2$ pour tout $(z,\zeta) \in \D(0,3 \eta/|c|)\times \overline{\D(0,3 \eta/|c|)}$.\\ 
\end{itemize}

\noindent Choisissons maintenant $\delta < \min\{A,A_{\eta},2 \eta/|c|\}$ et $\varepsilon_{1} \in ]0,\delta[$. Alors on a $\mathcal{L}_{\delta,\varepsilon_{1},\delta} \subset U_{\eta}$ (remarque \ref{rq34}) et l'inégalité (\ref{cond2}) ci-dessus est vérifiée.\\

Pour tout $|z| \leqslant \delta$ et $|\zeta|\in [\varepsilon_{1}, \delta]$ on a $(z, \zeta) \in K_{\delta,\varepsilon_{1}} \subset \overline{V}^{0}_{\eta,\delta}$. Soit maintenant $\ell \in \{1,...,k-1\}$ et $(z,\zeta)$ un point de $g^{\ell}(K)$, on a $|z| \geqslant \alpha$ et ceci entraîne que $(z,\zeta) \not \in \smash{\overline{V}}^{0}_{\eta,\delta}$. \\

En effet, supposons le contraire : la projection de $\smash{\overline{V}}^{0}_{\eta,\delta}$ sur la première coordonnée étant connexe, et comme $\delta < 2\eta/|c| < \alpha$, il devrait exister un élément $(z',\zeta') \in \smash{\overline{V}}^{0}_{\eta,\delta}$ avec $|z'| = 2\eta/|c|$, ce qui est impossible puisque dans ce cas $|q(z',\zeta')| > \displaystyle \frac{2\eta}{|c|}(|c|/2)=\eta$.

Ainsi, la composante connexe $\smash{\overline{V}}^{0}_{\eta,\delta}$ de $\overline{V}_{\eta,\delta}$ qui contient $K_{\delta,\varepsilon_{1}}$ ne rencontre aucune autre composante de $\mathcal{K}_{\delta,\varepsilon_{1}}$, ce qu'il fallait démontrer.\\

À partir de maintenant, on omet les indices $\delta$, $\varepsilon_{1}$.\\

\noindent \textbf{Étape 2 :} Montrons à présent que $\widehat{\mathcal{K}}^{0}=K$. Par l'étape 1, et comme $\widehat{\mathcal{K}}^{0}\subset \smash{\overline{V}}^{0}$, on sait que $\widehat{\mathcal{K}}^{0}$ ne rencontre pas d'autre composante de $\mathcal{K}$ que l'ensemble $K$ lui-même.\\
Soit $(z_{0},\zeta_{0}) \in \widehat{\mathcal{K}}^{0}\setminus K$. On suppose que $|\zeta_{0}|\geqslant \varepsilon_{1}$ (sinon $(z_{0},\zeta_{0}) \not \in \widehat{\mathcal{K}}$), donc nécessairement $|z_{0}|^{2}+|\zeta_{0}|^{2}>\delta^{2}$.  Comme la boule fermée $\overline{B}:=\overline{B(0,\delta)}$ est holomorphiquement convexe dans $\C^{2}$, il existe une fonction $h$ holomorphe sur $\C^{2}$ telle que $|h(z_{0},\zeta_{0})|> \displaystyle \|h\|_{\overline{B}}$. Notons respectivement $m_{0}$ et $m_{\overline{B}}$ les quantités $|h(z_{0},\zeta_{0})|$ et $\|h\|_{\overline{B}}$, ainsi que $m_{\widehat{\mathcal{K}}}$ la quantité $\|h\|_{\widehat{\mathcal{K}}}$, qui est finie puisque $\widehat{\mathcal{K}}$ est compact.\\

Considérons la fonction $\chi_{\widehat{\mathcal{K}}^{0}}$ définie sur $\widehat{\mathcal{K}}$ valant $1$ sur $\widehat{\mathcal{K}}^{0}$ (en particulier sur $K$) et $0$ sur $\widehat{\mathcal{K}}\setminus \widehat{\mathcal{K}}^{0}$ (en particulier sur $g^{\ell}(K)$ pour $\ell \not \equiv 0 [k]$). \\

Le théorème 6' p.~213 de \cite{Gunning:1965kx} nous dit que la fonction $\chi_{\widehat{\mathcal{K}}^{0}}$ est limite uniforme sur $\widehat{\mathcal{K}}$ de fonctions holomorphes sur $\C\times \Delta^{*}$. Soit donc $f$ une fonction holomorphe vérifiant $\|f-\chi_{\widehat{\mathcal{K}}^{0}}\|_{\widehat{\mathcal{K}}}< \varepsilon'$ avec $\varepsilon' < \displaystyle \min \left\{\frac{m_{0}}{m_{\widehat{\mathcal{K}}}+m_{0}}, \frac{m_{0}-m_{\overline{B}}}{m_{0}+m_{\overline{B}}}\right\}$ et appelons $F$ l'application $(z,\zeta) \mapsto h(z,\zeta) f(z,\zeta)$.\\

Pour $(z_{\ell},\zeta_{\ell}) \in g^{\ell}(K)$ (avec $\ell \in \{1,...,k-1\}$), on a l'inégalité \[|F(z_{\ell},\zeta_{\ell})|\leqslant \|F-h\chi_{\widehat{\mathcal{K}}^{0}} \|_{\widehat{\mathcal{K}}}+|h(z_{\ell},\zeta_{\ell})\chi_{\widehat{\mathcal{K}}^{0}}(z_{\ell},\zeta_{\ell}) |.\] Le second terme du membre de droite est nul ; quant au premier, il est majoré par $m_{\widehat{\mathcal{K}}}\varepsilon'$. De plus, on a $|F(z_{0},\zeta_{0})|=m_{0}|f(z_{0},\zeta_{0})|>m_{0}(1-\varepsilon'$) d'une part, et pour tout $(z,\zeta) \in K$ on a \[|F(z,\zeta)| \leqslant |h(z,\zeta)\left(f(z,\zeta)-\chi_{\widehat{\mathcal{K}}^{0}}(z,\zeta)\right)|+|h(z,\zeta)\chi_{\widehat{\mathcal{K}}^{0}}(z,\zeta)|\] donc $|F(z,\zeta)| \leqslant m_{\overline{B}} \varepsilon' + m_{\overline{B}}$ d'autre part. Le choix de $\varepsilon'$ nous assure que $\max\{m_{\widehat{\mathcal{K}}}\varepsilon', m_{\overline{B}}(\varepsilon' +1)\}<m_{0}(1-\varepsilon')$. Autrement dit, nous avons montré que $(z_{0},\zeta_{0}) \not \in \widehat{\mathcal{K}}$, d'où une contradiction. Ainsi, on a bien établi que $\widehat{\mathcal{K}}^{0}=K$. \\

\noindent \textbf{Étape 3 :} Il nous reste à conclure. Remarquons que l'enveloppe holomorphe convexe $\widehat{\mathcal{K}}$ de $\mathcal{K}$ est également stable par $g$ et supposons qu'il existe $\ell_{0} \in \{1,..., k-1\}$ et un point $(z_{\ell_{0}},\zeta_{\ell_{0}}) \in \widehat{\mathcal{K}}^{\ell_{0}}\setminus g^{\ell_{0}}(K)$. Alors on a les inclusions suivantes : \[K \subset g^{-\ell_{0}}(\widehat{\mathcal{K}}^{\ell_{0}}) \subset \widehat{\mathcal{K}}^{0},\] la dernière inclusion provenant du fait que la continuité de $g$ entraîne la connexité de $g^{-\ell_{0}}(\widehat{\mathcal{K}}^{\ell_{0}})$. On a donc $g^{-\ell_{0}}(z_{\ell_{0}},\zeta_{\ell_{0}}) \in \widehat{\mathcal{K}}^{0} = K$, d'où une contradiction.\\
Finalement, on a établi que $\displaystyle \bigcup_{\ell=0}^{k-1}\widehat{\mathcal{K}}^{\ell}=\mathcal{K}$. Comme $\mathcal{K}$ est une réunion de composantes connexes de $\widehat{\mathcal{K}}$, c'est un sous-ensemble ouvert et fermé de $\widehat{\mathcal{K}}$, donc holomorphiquement convexe par le corollaire 8 p.~214 de \cite{Gunning:1965kx}.}

Notons $\mathcal{O}(\C\times \Delta^{*})$ l'algèbre des fonctions holomorphes sur $\C\times \Delta^{*}$ et $\varphi^{*}(\mathcal{O}(\C\times \Delta^{*}))$ l'algèbre des éléments de $\mathcal{O}(\C\times \Delta^{*})$ invariants par le groupe $g^{\Z}$. Si $A$ est une algèbre de fonctions holomorphes, on note $\widehat{\mathcal{K}}^{A}$ l'enveloppe de $\mathcal{K}$ par rapport à l'algèbre $A$. On a montré que $\widehat{\mathcal{K}}^{\mathcal{O}(\C\times \Delta^{*})}=\mathcal{K}$.

\Corollaire{\label{corolEnv}On a $\widehat{\mathcal{K}}^{\varphi^{*}(\mathcal{O}(\C\times \Delta^{*}))}=\mathcal{K}$.}
\Preuve{En effet, pour $x \not \in \mathcal{K}$, on a : 
\begin{equation} \label{egEnv}
\widehat{(g^{\Z}.x)\cup \mathcal{K}}^{\mathcal{O}(\C\times\Delta^{*})}=(g^{\Z}.x) \cup \mathcal{K}.\end{equation}
Ceci découle du fait que si $p \not \in \mathcal{K}$, pour $q \not \in \{p\} \cup \mathcal{K}$, il existe $f_{1} \in \mathcal{O}(\C\times\Delta^{*})$ telle que $\|f_{1}\|_{\mathcal{K}}<f_{1}(q)$. Après avoir éventuellement multiplié $f_{1}$ par une constante, on peut supposer que $f_{1}(q)=1$. Comme $p\neq q$, il existe également une fonction $f_{2}\in \mathcal{O}(\C\times\Delta^{*})$ qui vérifie $f_{2}(p)=0$, $f_{2}(q) \neq 0$ et $\|f_{2}\|_{\mathcal{K}}\leqslant 1/2$ ; quitte à remplacer $f_{1}$ par des puissances d'elle-même, on peut supposer que $\|f_{1}\|_{\mathcal{K}}\leqslant |f_{2}(q)|$ et dans ce cas on a $\|f_{1}f_{2}\|_{\mathcal{K}\cup\{p\}} < |f_{1}(q)f_{2}(q)|$. Ainsi, on a $\widehat{\{p\}\cup \mathcal{K}}^{\mathcal{O}(\C\times\Delta^{*})}=\{p\} \cup \mathcal{K}$ ; par conséquent, en ajoutant un nombre fini de points à $\mathcal{K}$ l'ensemble obtenu reste holomorphiquement convexe, et on a bien l'égalité (\ref{egEnv}).

On considère alors la fonction $f$ qui vaut $1$ sur $g^{\Z}.x$ et $0$ sur $\mathcal{K}$, qui est holomorphe sur $(g^{\Z}.x) \cup \mathcal{K}$. Alors (théorème 6' p.~213 de \cite{Gunning:1965kx}) il existe une fonction $h \in \mathcal{O}(\C \times \Delta^*)$ telle que $\|f-h\|< 1/2$. 

En définissant la fonction holomorphe $\displaystyle H:=\frac{1}{k}\sum_{j=0}^{k-1}(h\circ g^{j})$, il sort que l'on a \[|H(x)-1|<1/2\] tandis que pour tout $y \in \mathcal{K}$, on a $|H(y)|<1/2$, donc $x \not \in \widehat{\mathcal{K}}^{\varphi^{*}(\mathcal{O}(\C\times \Delta^{*}))}$.}

\Corollaire{\label{corC}Soit $\delta$ un réel positif donné par la proposition \ref{holConv}. Alors, la paire \[\left(B(0,\delta)\setminus\{\zeta = 0\}, \varphi(B(0,\delta)\setminus\{\zeta = 0\})\right)\] est de Runge.}
\Preuve{On se donne un compact $A$ de $\varphi(B(0,\delta)\setminus\{\zeta = 0\})$, il est inclus dans un certain $\varphi(\mathcal{K}_{\delta-1/p,1/q})$ (pour $p$ et $q$ assez grands et avec $\delta > 1/p+1/q$). L'enveloppe de $A$ par rapport à l'algèbre des fonctions holomorphes sur $B(0,\delta)\setminus\{\zeta = 0\}$ est incluse dans $\varphi(\mathcal{K}_{\delta-1/p,1/q})$ par le corollaire \ref{corolEnv}, donc compacte. Ainsi $\varphi(B(0,\delta)\setminus\{\zeta = 0\})$ est holomorphiquement convexe par rapport aux fonctions holomorphes de $B(0,\delta)\setminus\{\zeta = 0\}$, ce qui nous donne la conclusion (\cite{Gunning:1965kx}, corollaire 9 p.~214).}

\SousChapitre{Une généralisation}
Soit $\varphi$ un germe de $(\C^{3},0)$ dans $(\C^{3},0)$ donné par \begin{equation}\label{formeNormaleGen}(z,\zeta,\xi) \mapsto (\lambda \zeta^r\xi^{s}z+P(\zeta,\xi),\zeta^{k},\xi^{\ell})\end{equation} où $\lambda \in \C^{*}$, $k, \ell, r, s \in \N$ avec $k,\ell>1$, $\pgcd(k,\ell)=1$ et $r,s>0$,  et \[P(\zeta,\xi)=\sum_{i_{1}=j_{1}}^{r}\sum_{i_{2}=j_{2}}^{s}c_{i_{1},i_{2}}\zeta^{i_{1}}\xi^{i_{2}}\] avec les conditions suivantes : $0<j_{1}<k$, $0<j_{2}<\ell$, $j_{1}\leqslant r$, $j_{2}\leqslant s$ et $c_{j_{1},j_{2}}\neq 0$.\\

Nous ajoutons une hypothèse supplémentaire, à savoir que pour tout $\varepsilon \in \U_{k}$ (racines $k$-ièmes de l'unité) et $\tau \in \U_{\ell}$ avec $\varepsilon \tau \neq 1$ \footnote{Comme $k$ et $\ell$ sont premiers entre eux, ceci revient à dire que $\varepsilon$ et $\tau$ ne sont pas simultanément égaux à $1$.}, il existe des entiers $n$ et $m$ et un polynôme $Q$ avec $Q(0,0)\neq0$, tels que l'on ait l'égalité : 
\begin{equation}\label{ok}
P(\zeta,\xi)-P(\varepsilon\zeta,\tau\xi)=\zeta^{n}\xi^{m}Q(\zeta,\xi).
\end{equation} 

Donnons quelques classes d'exemples de polynômes vérifiant cette dernière condition : 
\begin{itemize}
\item[1.] $P(\zeta,\xi) = \displaystyle \!\!\!\!\sum_{p=1}^{\min(r,s)}\!\!\!a_{p}\zeta^{p}\xi^{p}$ avec ou bien $\pgcd\{k,p~|~a_{p}\neq 0\} =1$, ou bien $\pgcd\{\ell,p~|~a_{p}\neq 0\} =1$,
\item[2.] $P(\zeta,\xi) = \displaystyle \zeta^{s'}\sum_{p=1}^{r}a_{p}\xi^{p} $ avec $\pgcd\{\ell,p~|~a_{p}\neq 0\} =1$ et $1\leqslant s' \leqslant s$,
\item[3.] $P$ de la forme précédente, mais en intervertissant les rôles de $\zeta$ et $\xi$.\\
\end{itemize}

Etant données $\varepsilon_{k}$ et $\tau_{\ell}$ deux racines primitives $k$-ième et $\ell$-ième de l'unité respectivement, notons $g$ l'automorphisme de $\C \times (\Delta^{*})^{2}$ qui à $(z,\zeta,\xi)$ associe \[(\underbrace{\varepsilon_{k}^{-r}\tau_{\ell}^{-s}z+\frac{P(\zeta,\xi)-P(\varepsilon_{k}\zeta,\tau_{\ell}\xi)}{\lambda\varepsilon_{k}^{r}\tau_{\ell}^{s}\zeta^{r}\xi^{s}}}_{\displaystyle a_{k,l}(z,\zeta,\xi)},\varepsilon_{k}\zeta, \tau_{\ell}\xi), \]
et $X$ l'ensemble $B(0,1)\setminus\{\zeta\xi=0\}$. La condition (\ref{ok}) permet d'adapter le raisonnement de la preuve de la proposition \ref{holConv} et de ses deux corollaires dans cette situation, en posant cette fois-ci $$q(z,\zeta,\xi)=z\displaystyle \prod_{i=1}^{k-1}\prod_{j=1}^{\ell-1}a_{k,\ell}(z,\zeta,\xi)\zeta^{n_{k}}\xi^{m_{\ell}}.$$

Ainsi la paire $(X,\varphi(X))$ est de Runge. On obtient alors une variété de Stein en recollant une infinité dénombrable de copies de $X \setminus \varphi(X)$ grâce à l'application $\varphi$. Il est possible de généraliser cette dernière construction en prenant un germe de $(\C^{n+1},0)$ dans lui-même, défini cette fois par $(z,\zeta_{1}, ..., \zeta_{n}) \mapsto (\lambda \zeta_{1}^{s_{1}}...\zeta_{n}^{s_{n}}z+P(\zeta_{1}, ..., \zeta_{n}),\zeta_{1}^{k_{1}}, ..., \zeta_{n}^{k_{n}})$ avec des conditions directement analogues à celles données ci-dessus.

\Chapitre{Invariants}
Revenons à présent à notre situation de départ. On note désormais $X$ la variété $\widetilde{S}\setminus \widetilde{D}$.
Etant donné un groupe $G$, on appelle espace $K(G,1)$ tout espace topologique connexe dont le groupe fondamental est isomorphe à $G$ et qui possède un revêtement universel contractile.

\Exemple{Le cercle unité $\S^{1}$ est un espace $K(\Z,1)$.}

Remarquons tout d'abord que la variété $X$ est un espace $K(\Z[\frac{1}{k}],1)$. En effet, $\pi_{1}(X) \cong \Z[\frac{1}{k}]$ et son revêtement universel $\C\times\H$ (c.f. \cite{Dloussky:1999aa} et \cite{Favre:2002aa}) est contractile.

Le théorème I de \cite{Eilenberg:1945aa} (pp.~482-483) nous dit alors que les groupes de cohomologie de $X$ sont isomorphes à ceux du groupe $\Z[\frac{1}{k}]$, c'est-à-dire que pour tout $n \in \N$ et pour tout groupe $G$, on a un isomorphisme entre $H^{n}(X,G)$ et $H^{n}(\Z[\frac{1}{k}],G)$. De plus, on sait (loc. cit. pp.~488-489) que le groupe $H^{2}(\Z[\frac{1}{k}],G)$ est isomorphe au groupe des extensions centrales de $\Z[\frac{1}{k}]$ par $G$. Une extension centrale est la donnée d'une extension de groupe \[0 \rightarrow G \overset{i}{\rightarrow} E \overset{p}{\rightarrow} \textstyle{\Z[\frac{1}{k}]} \rightarrow 0\]
où $E$ est un groupe avec $i(G) \subset Z(E)$, le centre de $E$.\\

Nous sommes maintenant en mesure de prouver la 
\Proposition{Le groupe $H^{2}(X,\C)$ est trivial.}
\Preuve{Par ce qui précède, il suffit de montrer qu'une extension centrale $E$ de $\Z[\frac{1}{k}]$ par $\C$ est nécessairement triviale, i.e. isomorphe au produit cartésien $\C\times \Z[\frac{1}{k}]$. Soit donc $E$ une telle extension :  \[0 \rightarrow \C \overset{i}{\rightarrow} E \overset{p}{\rightarrow} \textstyle{\Z[\frac{1}{k}]} \rightarrow 0.\] Montrons que $E$ est abélien. Soient $x, y \in E$ et $a \in \N$ tels que $p(x)$ et $p(y)$ appartiennent tous deux à $\frac{1}{k^{a}}\Z := \{\frac{n}{k^{a}}, n \in \Z\}$ qui est un sous-groupe de $\Z[\frac{1}{k}]$ isomorphe à $\Z$. \\
L'extension $E$ induit une extension $F:=p^{-1}(\frac{1}{k^{a}}\Z)$ de $\frac{1}{k^{a}}\Z$ par $\C$, donc une extension de $\Z$ par $\C$ :
 \[0 \rightarrow \C \overset{i'}{\rightarrow} F \overset{p'}{\rightarrow} \textstyle{\Z} \rightarrow 0\]
 
Il existe une section $s : \Z \rightarrow F$ (on choisit $s(1) \in p'^{-1}(1)$ et on pose $s(n)=ns(1)$ pour $n\in\Z$) donc $F$ est produit semi-direct de $\Z$ par $\C$, donné par $\sigma \in \Hom(\Z,\Aut(\C))$. L'extension $F$ étant elle aussi centrale, $\sigma \equiv 1$ est l'unique possibilité, i.e. $F$ est abélien (il est isomorphe à $\C\times\Z$) donc $x$ et $y$ commutent. Ainsi, $E$ est abélien. 

Il existe des sections $s : \Z[\frac{1}{k}]\mapsto E$. Pour construire l'une d'elles, fixons $x_{0}\in p^{-1}(1)$. Comme $\Z[\frac{1}{k}] \cong E/i(\C)$, il existe $x'_{1}\in p^{-1}(1/k)$ tel que $k x'_{1} = x_{0}+i(w)$ avec $w\in \C$. On pose alors $x_{1}:=x'_{1}-i(w/k)$ et on a $kx_{1}=x_{0}$ ; on définit ainsi par récurrence les $x_{i} \in p^{-1}(1/k^{i})$ vérifiant $kx_{i+1}=x_{i}$, et notre section est donnée par $s(n/k^{a})=nx_{a}$ pour $n \in \Z$ et $a \in \N$. L'existence d'une telle section nous dit que $E$ est isomorphe au produit semi-direct $\C\rtimes\Z[\frac{1}{k}]$ donné par $\sigma \in \Hom(\Z[\frac{1}{k}],\Aut(\C))$. Le groupe $E$ étant abélien, on a nécessairement $\sigma \equiv 1$, i.e. $E$ est isomorphe au produit $\C\times \Z[\frac{1}{k}]$.}
\Remarque{Le groupe $H^{2}(X,\Z)$ n'est pas trivial ; il contient des éléments de torsion et des éléments qui ne sont pas d'ordre fini. Le morphisme $\rho$ de la preuve du lemme \ref{lem2} ci-après fournit un exemple d'élément de torsion, puisqu'on peut voir que $\rho^{k-1}$ admet un logarithme (voir définition \ref{log} ci-dessous). Pour ce qui est des éléments qui ne sont pas d'ordre fini, donnons-en un exemple. Considérons le groupe $\Z[\frac{1}{6}]$. On a un isomorphisme de groupes $\varphi : \Z[\frac{1}{6}]/\Z[\frac{1}{2}]\cong \Z[\frac{1}{3}]/\Z$ et une injection $i$ de ce groupe (le $3$-groupe de Prüfer) dans $\S^{1}$. Alors on peut voir que $\rho:=i\circ \varphi$ n'est pas d'ordre fini dans $H^{2}(X,\Z)$.
Ainsi, le groupe $H^{2}(X,\Z)$ possède des éléments d'ordre infini dont l'image est nulle dans $H^{2}(X,\Q)$, ceci est conséquence du fait que le groupe $H_{1}(X,\Z) \cong \Z[\frac{1}{k}]$ n'est pas finiment engendré (voir \cite{Bott:1972aa}, théorème 4 p.~144).}

Étant donnée une surface intermédiaire $S$ et son germe associé sous la forme normale (\ref{formeNormale}), on définit l'indice de $S$ comme le plus petit entier $m$ tel que $k-1$ divise $m s$ (voir \cite{Oeljeklaus:2008aa}). 

Il existe un feuilletage holomorphe $\mathcal{F}$ sur $X$ défini par la $1$-forme holomorphe  $\omega=\displaystyle \frac{d\zeta}{\zeta}$, qui ne s'annule nulle part (c.f. \cite{Dloussky:2001aa}). De façon équivalente, les feuilles de ce feuilletage sont les ensembles $\{\zeta=const.\}$.
Dans le cas où $S$ est d'indice $1$, i.e. lorsque $k-1$ divise $s$, il existe un champ de vecteurs tangent à ce feuilletage qui ne s'annule nulle part, autrement dit on a le 
\Lemme{\label{lem0}Lorsque la surface $S$ est d'indice $1$, le fibré tangent au feuilletage $T\mathcal{F}$ est holomorphiquement trivialisable.}
\Preuve{Pour prouver cela, nous montrons qu'il suffit de considérer le champ de vecteurs $V$ sur $X$ induit par le champ de vecteurs $\widetilde{V}=\zeta^{\frac{s}{k-1}}\frac{\partial}{\partial z}$ sur $\C\times\Delta^{*}$, tangent au feuilletage de $\C\times\Delta^{*}$ défini par $\omega$.\\

En effet, d'une part on remarque que $X$ est le quotient de $\C\times \Delta^{*}$ par $G$ où $G \cong \Z[\frac{1}{k}]/\Z$ est le groupe formé des automorphismes de $\C\times \Delta^{*}$ de la forme $g_{k^{n}}^{\ell}(z,\zeta)=(z\varepsilon_{k^{n}}^{-\ell s \frac{k^{n}-1}{k-1}}+$\[\frac{\displaystyle \sum_{i=0}^{n-1}\lambda^{n-i-1}\zeta^{sk^{i+1}\frac{k^{n-i-1}-1}{k-1}}(P(\zeta^{k^{i}})\!-\!\varepsilon_{k^{n}}^{-\ell s k^{i+1}\frac{k^{n-i-1}-1}{k-1}}P((\varepsilon_{k^{n}}^{\ell}\zeta)^{k^{i}}))}{\lambda^{n}(\varepsilon_{k^{n}}^{\ell}\zeta)^{s\frac{k^{n}-1}{k-1}}},\varepsilon_{k^{n}}^{\ell}\zeta)\]
pour $n \in \N, \ell \in \{0,..., k^{n}-1\}$ et avec $\varepsilon_{k^{n}}=e^{\frac{2i\pi}{k^{n}}}$. Ceci provient du fait que $X$ est le quotient de $\C\times \H_{g}$ par le groupe $\{\gamma^{n}\gamma_{1}^{\ell}\gamma^{-n}~|~n,\ell \in \Z\}\cong\Z[\frac{1}{k}]$ où $\H_{g}=\{w \in \C~|~\Re(w)<0\}$, $\gamma(z,w)=(\lambda z e^{sw}+P(e^{w}) ,kw)$ et $\gamma_{1}(z,w)=(z,w+2i\pi)$ (voir \cite{Dloussky:2001aa}, proposition 2.3 et section 4). On considère alors le quotient par le sous-groupe $\{\gamma^{n}\gamma_{1}^{k^{n}\ell}\gamma^{-n}~|~n,\ell \in \Z\}\cong\Z$ ce qui nous donne bien $X = (\C\times\Delta^{*})/G$.\\

D'autre part, un champ de vecteurs $\widetilde{V}$ défini sur $\C\times\Delta^{*}$ induit un champ de vecteurs tangent à $X$ lorsqu'il est invariant par le groupe $G$, i.e. s'il vérifie : 
\begin{equation}\label{condVcTg}
D(g_{k^{n}}^{\ell})_{_{z,\zeta}}(\widetilde{V}(z,\zeta))=\widetilde{V}(g^{\ell}_{k^{n}}(z,\zeta)).
\end{equation}
Cette condition est bien vérifiée par $\widetilde{V}$ puisque l'on a l'égalité $\varepsilon_{k^{n}}^{-\ell s \frac{k^{n}-1}{k-1}}\zeta^{\frac{s}{k-1}}=(\varepsilon_{k^{n}}^{\ell}\zeta)^{\frac{s}{k-1}}$. Comme $\omega(\widetilde{V})=0$, on a bien montré que $V\in H^{0}(X,T\mathcal{F})$.}


\Remarque{On peut montrer qu'un champ de vecteurs sur $X$ de la forme $f(z,\zeta)\frac{\partial}{\partial z}$ (où $f$ est une holomorphe ne s'annulant nulle part) existe bien si et seulement si la surface $S$ est  d'indice $1$, autrement dit on a une équivalence dans le lemme précédent. C'est une conséquence de la condition (\ref{condVcTg}) et le raisonnement est analogue à celui qui sera fait dans le lemme \ref{lem2}.}

La trivialité du fibré $T\mathcal{F}$ entraine celle du fibré tangent $TX$, ce que nous voyons à présent.

\Lemme{\label{lem1}Lorsque le fibré $T\mathcal{F}$ est holomorphiquement trivialisable, le fibré tangent holomorphe $TX$ de $X$ l'est aussi.}
\Preuve{Étant donné que nous avons une section holomorphe globale $V$ de $T\mathcal{F}$, il nous suffit d'exhiber un deuxième champ de vecteurs global, linéairement indépendant de $V$ en chaque point. Par définition, on peut trouver un recouvrement de $X$ par des ouverts $U_{i}$ et sur chacun d'eux un champ de vecteurs $W_{i}$ qui soit linéairement indépendant de $V$ sur $U_{i}$. Quitte à remplacer $W_{i}$ par $W_{i}/\omega(W_{i})$ on peut supposer que $\omega(W_{i})\equiv 1$ sur $U_{i}$, de sorte que $\omega(W_{i,j})=0$ sur $U_{i,j}:=U_{i}\cap U_{j}$, où l'on a posé $W_{i,j}:=W_{i}-W_{j}$. La famille $(W_{i,j})$ forme donc un cocyle de $H^{1}(X,T\mathcal{F})$ qui est aussi un cobord par le théorème B de Cartan. Ainsi il existe un champ de vecteurs $Z_{i}$ sur chaque $U_{i}$ tel que $Z_{i}-Z_{j}=W_{i,j}$. Posons $\widetilde{Y}_{i}:=W_{i}-Z_{i}$, de sorte que  $\widetilde{Y}_{i}=\widetilde{Y}_{j}$ sur $U_{i,j}$, i.e. les $\widetilde{Y}_{i}$ se recollent en une section holomorphe globale de $TX$. Nous avons deux champs de vecteurs $V$ et $\widetilde{Y}$ vérifiant $\omega(V)\equiv 0$ et $\omega(\widetilde{Y})\equiv 1$, ce qui nous assure qu'ils sont linéairement indépendants en chaque point de la variété étudiée.}

Nous voulons à présent établir un lien entre le fait que $S$ soit d'indice $1$ et la trivialité du fibré canonique de $X$.

On considère la suite exacte courte $0 \rightarrow \Z \rightarrow \C \rightarrow \C^{*} \rightarrow 0$ de faisceaux, qui donne lieu à la suite exacte longue de cohomologie $ \cdots \rightarrow H^{1}(X,\Z) \rightarrow H^{1}(X,\C) \rightarrow H^{1}(X,\C^{*}) \rightarrow H^{2}(X,\Z) \rightarrow H^{2}(X,\C) \rightarrow \cdots$. Toujours d'après le théorème I de \cite{Eilenberg:1945aa}, on a les isomorphismes $$\textstyle H^{1}(X,\Z) \cong \Hom(\Z\!\left[\frac{1}{k}\right],\Z)=0,$$ $$\textstyle H^{1}(X,\C^{*}) \cong \Hom(\Z\!\left[\frac{1}{k}\right],\C^{*})=\C^{*}$$ et $$\textstyle H^{1}(X,\C) \cong \Hom(\Z\!\left[\frac{1}{k}\right],\C)\cong \C.$$ D'autre part le groupe $H^{2}(X,\C)$ est trivial, d'où l'on tire finalement la suite exacte courte $$0 \rightarrow H^{1}(X,\C) \overset{e^{2i\pi \cdot}}{\rightarrow} H^{1}(X,\C^{*}) \overset{c}{\rightarrow} H^{2}(X,\Z) \rightarrow 0.$$

\Definition{\label{log}On dira qu'un élément $\rho$ de $H^{1}(X,\C^{*})$ admet un logarithme lorsqu'il existe un morphisme $\rho'$ de $\Z[\frac{1}{k}]$ dans $\C$ tel que $e^{2i\pi\rho'}=\rho$.} Ainsi, l'image par $c$ d'un élément $\rho \in H^{1}(X,\C^{*})$ est triviale dans $H^{2}(X,\Z)$ si et seulement si $\rho$ admet un logarithme $\rho'$.

\Lemme{\label{lem2}Si le fibré canonique de $X$ est holomorphiquement trivialisable, alors la surface $S$ est d'indice $1$.}
\Preuve{Le fibré canonique de $X$ est le fibré des $2$-formes holomorphes sur $X$ ; raisonnons par l'absurde et supposons qu'il est holomorphiquement trivialisable et que la surface $S$ n'est pas d'indice $1$. Alors il existe une $2$-forme holomorphe globale sur $X$ qui ne s'annule nulle part. 
Une telle forme provient d'une $2$-forme holomorphe sur le revêtement $\C\times \Delta^{*}$ de $X$ donnée par $f(z,\zeta)dz\wedge d\zeta$ (où $f$ est une fonction holomorphe sur $\C\times \Delta^{*}$ qui ne s'annule nulle part) qui soit stable par le groupe $G=\{g_{k^{n}}^{\alpha}~|~ n\in \N, \alpha \in \{0,...,k^{n}-1\} \}$, i.e. vérifie l'équation 
\[(g_{k^{n}}^{\alpha})^{*}(f(z,\zeta)dz\wedge d\zeta)=f(z,\zeta)dz\wedge d\zeta\] (pour tout $n\in\N$ et $\alpha\in\{0,...,k^{n}-1\}$). Ceci donne la condition suivante sur la fonction $f$ : 
\begin{equation}\label{condComp}e^{2i\pi\frac{\alpha}{k^{n}} (s\frac{k^{n}-1}{k-1}-1)}f(z,\zeta) = f(g_{k^{n}}^{\alpha}(z,\zeta)).
\end{equation}

Considérons l'homomorphisme de groupes \[\application{\rho}{\Z[\frac{1}{k}]}{\S^{1}.}{\frac{\alpha}{k^{n}}}{e^{2i\pi\frac{\alpha}{k^{n}} (s\frac{k^{n}-1}{k-1}-1)}}\] Il induit un fibré plat $L_{\rho}$ au-dessus de $X$, qui est holomorphiquement trivialisable si et seulement si $\rho$ admet un logarithme, puisque $H^{1}(X,\mathcal{O}^{*})\cong H^{2}(X,\Z)$ car $X$ est de Stein. Étant donné que la fonction $f$ vérifie la condition (\ref{condComp}) ci-dessus, elle définit une section holomorphe du fibré plat $L_{\rho}$ au-dessus de $X$. 

Ainsi pour pouvoir aboutir à une contradiction, il nous reste à voir que $\rho$ n'admet pas de logarithme (et donc qu'une telle fonction $f$ n'existe pas).

Remarquons tout d'abord que l'application $\sigma : \frac{\alpha}{k^{n}}\mapsto e^{2i\pi\frac{-\alpha}{k^{n}}(\frac{s}{k-1}+1)}$ est un homomorphisme de $\Z[\frac{1}{k}]$ dans $\S^{1}$ qui admet un logarithme. Ainsi, $\rho$ admet un logarithme si et seulement si l'homomorphisme $\varphi := \rho/\sigma :  \frac{\alpha}{k^{n}} \mapsto e^{2i\pi \frac{\alpha s}{k-1}}$ admet un logarithme.

Soit $m$ l'indice de la surface $S$. Comme $k-1$ n'est pas un diviseur de $s$, le noyau de $\varphi$ est précisément $m\Z[\frac{1}{k}]$ et cet homomorphisme n'admet donc pas de logarithme. En effet, si un tel morphisme $\rho'$ existait, sa restriction à $m\Z[\frac{1}{k}]$ serait un homomorphisme à valeurs dans $\Z$, nécessairement trivial. On aurait alors $m.\rho\left(\frac{1}{k^{n}}\right)=0$, i.e. $\rho\left(\frac{1}{k^{n}}\right)=0$ pour tout $n\in \N$.}

\Remarque{On a en fait une équivalence dans le lemme précédent. Lorsque la surface $S$ est d'indice $1$, on considère la forme $\zeta^{-(\frac{s}{k-1}+1)}dz \wedge d\zeta$, qui trivialise le fibré canonique.}

Les trois lemmes précédents ont en particulier comme conséquence la 
\Proposition{\label{indTriv}Soient $S$ une surface intermédiaire et $X=\widetilde{S}\setminus\widetilde{D}$. Les trois assertions suivantes sont équivalentes :
\begin{itemize}
\item[1.] La surface $S$ est d'indice $1$,
\item[2.] Le fibré tangent au feuilletage $T\mathcal{F}$ de $X$ est holomorphiquement trivialisable,
\item[3.] Le fibré tangent holomorphe $TX$ de $X$ est holomorphiquement trivialisable.
\end{itemize}
}
\Preuve{Vu les lemmes \ref{lem0} et \ref{lem1}, il suffit de montrer que la troisième assertion entraine la première. C'est une conséquence du lemme \ref{lem2}, car si $S$ n'est pas d'indice $1$, le fibré canonique de $X$ n'est pas holomorphiquement trivialisable. Dans ce cas, le fibré cotangent de $X$ et donc le fibré tangent $TX$ ne le sont pas non plus.}

\Remarque{Le problème suivant demeure non résolu actuellement (voir~\cite{Forstneric:2003aa}) : une variété de Stein de dimension $n$ dont le fibré tangent holomorphe est holomorphiquement trivialisable est-elle nécessairement un domaine de Riemann au-dessus de $\C^{n}$ ? Nous ne connaissons pas la réponse pour les surfaces de Stein que l'on vient de considérer.}

\bibliographystyle{amsplain}
\bibliography{biblio}

\end{document}